\documentclass[11pt]{article}
\title{Curvature Estimates in Asymptotically Flat Manifolds\\
of Positive Scalar Curvature}
\author{Felix Finster and Ines Kath}
\date{January 2001}

\newtheorem{Def}{Definition}[section]
\newtheorem{Thm}[Def]{Theorem}

\newtheorem{Lemma}[Def]{Lemma}

\newtheorem{Corollary}[Def]{Corollary}
\newcommand{\Proof}{{\em{Proof: }}}
\newcommand{\QED}{\ \hfill $\FBox$ \\[1em]}
\newcommand{\spc}{\;\;\;\;\;\;\;\;\;\;}

\newfont{\fra}{eufm10 scaled 1095}
\newfont{\Bb}{msbm10 scaled 1095}
\newcommand\C{{\mbox{\Bb C}}}
\newcommand\R{{\mbox{\Bb R}}}

\newcommand{\1}{\mbox{\rm 1 \hspace{-1.05 em} 1}}

\newcommand{\Tr}{\mbox{Tr\/}}
\newcommand{\Span}{\mathop{{\rm span}}}
\newcommand{\Ric}{\mathop{{\rm Ric}}}
\newcommand{\Div}{\mathop{{\rm div}}}
\newcommand{\Grad}{\mathop{{\rm grad}}}
\newcommand{\FBox}{\rule{2mm}{2.25mm}}

\newcommand\ph{{\varphi}}

\begin{document}

\maketitle
\begin{abstract}
We consider an asymptotically flat Riemannian spin manifold of positive scalar
curvature.  An inequality is derived which bounds the Riemann tensor in terms of
the total mass and quantifies in which sense curvature must become small when
the total mass tends to zero.
\end{abstract}

\section{Introduction}
Suppose that $(M^n, g)$ is an asymptotically flat Riemannian spin manifold
of positive scalar curvature. The positive mass theorem~\cite{SY, W, Br}
states that the total mass of the manifold is always positive, and is
zero if and only if the
manifold is flat. This result suggests that there should be an inequality
which bounds the Riemann tensor in terms of the total mass and implies
that curvature must become small when the total mass tends to zero.
In~\cite{BF} such curvature estimates were derived in the context of General
Relativity for $3$-manifolds being hypersurfaces in a Lorentzian manifold.
In the present paper, we study the problem more generally on a
Riemannian manifold of dimension $n \geq 3$. Our curvature estimates then
give a quantitative relation between the local geometry and global
properties of the manifold.

The main difficulty in higher dimensions is to bound the Weyl tensor (which
for $n=3$ vanishes identically). Our basic strategy for controlling the Weyl
tensor can be understood from the following simple consideration.
The existence of a parallel spinor in an open set $U \subset M$ implies
that the manifold is Ricci flat in $U$. Thus it is reasonable that by
getting suitable estimates for the derivatives of a spinor, one can bound
all components of the Ricci tensor. This method is used in~\cite{BF}, where
a solution of the Dirac equation is analyzed using the Weitzenb\"ock formula.
But the local existence of a parallel spinor does not imply that
the Weyl tensor vanishes. This is the underlying reason why in dimension
$n>3$, our estimates cannot be obtained by looking at one spinor, but
we must consider a family $(\psi^i)_{i=1,\ldots,2^{[n/2]}}$ of solutions
of the Dirac equation.
Out of these solutions we form the so-called spinor operator $P_x$.  The
curvature tensor can be bounded in terms of suitable derivatives of $P_x$, and
an integration-by-parts argument, the Weitzenb\"ock formula, and a-priori
estimates for the spinor operator give the desired result.

We now give the precise statement of our result.  For simplicity, we consider
only one asymptotically flat end.  The following definition immediately
generalizes that used in~\cite{BF}; for a slightly more general definition
see~\cite{Br}.
\begin{Def}\label{def1}
A Riemannian manifold $(M^n, g)$, $n \geq 3$,
is said to be {\em{asymptotically flat}} if
there is a compact set $K \subset M$ and a diffeomorphism $\Phi$ which maps
$M \setminus K$ to the region $\R^n \setminus B_{r_0}(0)$ outside a ball of
radius $r_0$. Under this diffeomorphism, the metric should be of the form
\[ (\Phi_* g)_{ij} \;=\; \delta_{ij} \:+\: {\cal{O}}(r^{2-n}) \;\;\;,\;\;\;\;\;
\partial_k (\Phi_* g)_{ij} \;=\; {\cal{O}}(r^{1-n}) \;\;\;,\;\;\;\;\;
\partial_{kl} (\Phi_* g)_{ij} \;=\; {\cal{O}}(r^{-n}) \;. \]
Furthermore, scalar curvature should be in $L^1(M)$.
\end{Def}
For an asymptotically flat manifold the total mass $m$ is defined by
\begin{equation}
m \;=\; \frac{1}{c(n)} \lim_{\rho \rightarrow \infty}
\int_{S_\rho} (\partial_j (\Phi_* g)_{ij} - \partial_i (\Phi_* g)_{jj})
\:d\Omega^i \;, \label{eq:A}
\end{equation}
where $c(n)>0$ is a normalization constant (which can be chosen arbitrarily),
$d\Omega$ is the volume form on the sphere $S_\rho \subset \R^n$ of radius
$\rho$, and $d\Omega^i$ denotes the product of $d\Omega$ by the i-th
component of the normal vector of $S_\rho$. As shown in~\cite{Br}, this
definition is independent of the
choice of $\Phi$. For our estimates, we also need the isoperimetric
constant $k$ given by
\[ k \;=\; \inf \frac{A}{V^\frac{n-1}{n}} \;, \]
where the infimum is taken over all smooth regions $D$, $V$ is the $n$-volume
of $D$, and $A$ is the $(n-1)$-volume of the boundary of $D$.

\begin{Thm} \label{thm1}
Let $M^n$ be an asymptotically flat Riemannian spin manifold of
positive scalar curvature. Then there are positive constants $c_1$, $c_2$,
and $c_3$ depending only on $n$ as well as a set $D$ with
\begin{equation}
    \mu(D) \;\leq\; \left( \frac{c_3\:m}{k^2} \right)^{\frac{n}{n-2}}
    \label{eq:ta}
\end{equation}
such that for all positive $\eta \in C^\infty(M)$ with $\sup_M(|\eta|+|\nabla \eta|)
<\infty$ the following inequality holds,
\[ \int_{M \setminus D} \eta\: |R|^2 \:d\mu \;\leq\;
m \:c_1\: \sup_M(|\eta R| + |\Delta \eta|) \:+\:
\sqrt{m}\:c_2\: \|\eta \:\nabla R\|_2 \;. \]
Here $R$ is the Riemann tensor, $m$ the total mass, and $k$ the
isoperimetric constant, $\|\cdot\|_2$ denotes the $L^2$-Norm on $M$.
\end{Thm}

One application of the above theorem is to a continuous family
$(M_\lambda)_{\lambda \geq 0}$ of manifolds (e.g.\ obtained by a flow of the
metric). If $\lim_{\lambda \rightarrow \infty} m(\lambda)=0$,
$\sup_M |R|$ and $\|\nabla R \|_2$ are uniformly bounded, and the
isoperimetric constant is bounded away from zero, then our theorem
implies that the Riemann tensor converges to zero pointwise almost everywhere,
and thus the manifold becomes flat. For other applications see~\cite{BF}.

\section{Spinors, the Dirac Equation}
In this section we recall some basic facts concerning spinors and the
Dirac equation; for details see \cite{B, LM, Br}.
Let $(M^n,g)$ be an $n$-dimensional (oriented) Riemannian spin
manifold with spin structure $Q$ and spinor bundle
\[S=Q\times_{Spin(n)}\Delta_n \;,\]
which is associated with $Q$ by the spinor
representation $\Delta_n$.  As a vector space, $\Delta_n$ is equal to $\C^N$,
$N:=2^{[n/2]}$.  The canonical Hermitian product on $\Delta_n$ defines a
complex scalar product $\langle\cdot \,,\cdot \rangle$ on $S$.
We denote the real part of this scalar product by $(\cdot \,,\cdot )$.
The Levi-Civita connection $\nabla$ on $(M^n,g)$ induces a covariant
derivative in $S$ called spinor derivative, which we denote again by
$\nabla$.  The spinor derivative is compatible with
$\langle\cdot \,,\cdot \rangle$, i.e.
\[X\langle \ph,\psi\rangle=\langle \nabla_X\ph,\psi\rangle + \langle
\ph,\nabla_X\psi\rangle\]
for all sections $\ph, \psi$ in $S$ and all vector
fields $X$ on $M$.  Its curvature tensor
\begin{equation}
R^S(X,Y) = \nabla_X \nabla_Y -\nabla_Y\nabla_X - \nabla_{[X,Y]}
\label{eq:2Z}
\end{equation}
is given locally by
\begin{equation}
R^S(X,Y)\: \psi=\frac{1}{4} \sum_{\alpha,\beta=1}^n R(X,Y,s_\alpha,s_\beta)
s_\alpha \cdot s_\beta \cdot \psi, \label{eq:2Y}
\end{equation}
where $s_1,\dots ,s_n$ is a local orthonormal frame on $(M^n,g)$, $R$ is the
Riemannian curvature tensor of $(M^n,g)$ and ``$\cdot$'' denotes Clifford
multiplication of a spinor by a vector. The Clifford multiplication satisfies
the anti-commutation relations
\begin{equation}\label{eq:acomm}
X\cdot Y\cdot\psi +Y\cdot X\cdot\psi = -2g(X,Y)\: \psi
\end{equation}
and is anti-Hermitian,
\[\langle X \cdot \ph,\psi\rangle + \langle \ph,X \cdot \psi\rangle =0.\]

The Dirac operator $D$ on $(M^n,g)$ is the composition of
the spinor derivative $\nabla: \Gamma(S)\longrightarrow \Gamma(TM\otimes S)$
and the Clifford multiplication, i.e.\ locally
\[ D \;=\; \sum_{\alpha=1}^n s_\alpha \cdot\nabla_{s_\alpha} \;. \]
The square $D^2$ of the Dirac operator
satisfies the Weitzenb\"ock formula
\begin{equation} D^2=\Delta^S + \frac{\tau}{4} \;, \label{wb} \end{equation}
where $\Delta^S=-\sum_{\alpha=1}^n (\nabla_{s_\alpha} \nabla_{s_\alpha}+
(\Div {s_\alpha}) \nabla_{s_\alpha})$
denotes the Bochner-Laplace operator with respect to the
spinor derivative and $\tau$ is the scalar curvature of $(M^n,g)$.

Now assume furthermore that $M$ is asymptotically flat
and has positive scalar curvature. In the coordinates induced by
the diffeomorphism $\Phi$ of Definition~\ref{def1}, we choose a constant spinor
$\psi_0$ and consider the boundary problem
\begin{equation}
D \psi \;=\; 0 \;\;\;,\spc \lim_{|x| \rightarrow \infty} \psi(x) \;=\;
\psi_0 \;. \label{eq:1C}
\end{equation}
It is shown in~\cite{Br} that this boundary problem has a unique weak
solution $\psi$. Using the asymptotic form of the metric in
Definition~\ref{def1}
and elliptic regularity theory, it follows that $\psi$ is even smooth and
decays at infinity like
\begin{equation}\label{psi}
 \psi \;=\; \psi_0 \:+\: {\cal{O}}(r^{2-n}) \;\;\;,\;\;\;\;\;
\partial_k \psi \;=\; {\cal{O}}(r^{1-n}) \;\;\;,\;\;\;\;\;
\partial_{kl}\psi \;=\; {\cal{O}}(r^{-n}) \;. \end{equation}
This solution of the Dirac equation can be used to prove the positive mass
theorem~\cite{W, Br}, as we now briefly outline. Consider the vector field
\begin{equation}
X \;=\; \frac{1}{2} \: \Grad |\psi|^2 \;. \label{eq:27a}
\end{equation}
Using the Weitzenb\"ock formula and the positivity of scalar curvature, we can
estimate the divergence of $X$ as follows,
\begin{eqnarray}
\Div X &=& \sum_{\alpha=1}^n s_\alpha \:(\nabla_{s_\alpha}
\psi, \:\psi) \;=\; -(\Delta^S \psi,\: \psi) \:+\: |\nabla \psi|^2 \nonumber \\
&=& \frac{\tau}{4} \:|\psi|^2 \:+\: |\nabla \psi|^2 \;\geq\;
|\nabla \psi|^2 \;. \label{eq:27b}
\end{eqnarray}
We introduce the balls $B_\rho$ and spheres $S_\rho$ of radius $\rho>r_0$ by
\begin{equation}
    D_\rho \;=\; K\cup \Phi^{-1}(B_\rho(0)\setminus B_{r_0}(0))\;\;\;,\spc
    S_\rho \;=\; \partial D_\rho
    \label{eq:Dr}
\end{equation}
with $\Phi$, $r_0$, and $K$ as in Definition~\ref{def1}, $\rho>r_0$, and
denote the outer normal on $S_\rho$ by $\nu$.
We integrate over $D_\rho$ and apply Gauss' theorem,
\begin{eqnarray}
\int_{S_\rho} g( X, \nu) \:dS_\rho \;=\;
\int_{D_\rho} \Div X \:d\mu \;\geq\;
\int_{D_\rho} |\nabla \psi|^2 \:d\mu \;. \label{eq:1B}
\end{eqnarray}
An asymptotic expansion near infinity~\cite{Br} shows that as $\rho \rightarrow
\infty$, the left side of~(\ref{eq:1B}) can up to a constant be identified
with the boundary integral in~(\ref{eq:A}); more precisely,
\[ \lim_{\rho \rightarrow \infty} 4
\int_{S_\rho} g( X, \nu) \:dS\rho
\;=\; c(n) \: |\psi_0|^2 \:m \;. \]
Hence in~(\ref{eq:1B}) we can take the limit $\rho \rightarrow \infty$
to obtain
\begin{eqnarray}
c(n)\: |\psi_0|^2 \:m \;\geq\; 4\:\|\nabla \psi\|_2^2 \;.
\label{eq:1C2}
\end{eqnarray}
This inequality shows that $m \geq 0$. If $m=0$, (\ref{eq:1C2}) yields that
for any $\psi_0$, there is a parallel spinor with $\lim_{|x| \rightarrow
\infty} \psi(x) = \psi_0$, and this implies that the manifold is flat.

The inequality~(\ref{eq:27b}) immediately gives an a-priori bound for the
spinor, which will be very useful later on. Namely, (\ref{eq:27a}) and
(\ref{eq:27b}) imply that
\[ \Delta |\psi|^2 \;=\; -2 \:\Div X \;\leq\; 0\;. \]
Thus $|\psi|^2$ is sub-harmonic, and the maximum principle yields that for
every solution of the boundary problem~(\ref{eq:1C}),
\begin{equation}
\sup_M|\psi|\le |\psi_0|\,.
\label{eq:1E}
\end{equation}

\section{The Integration-by-Parts Argument}
In this section we derive an $L^2$-bound for the second derivative of a solution
of the Dirac equation~(\ref{eq:1C}). The argument is similar to that in
dimension three~\cite{BF}. We give it in some detail, using the formalism of
orthonormal frames~\cite{B}.

We define a vector field $Y$ on $M$ by
\[Y=\frac 12 \Grad |\nabla \psi |^2\,.\]
\begin{Lemma} \label{L41}
For any local orthonormal frame $s_1,\dots,s_n$ on $M$,
\begin{eqnarray*} \Div Y &=& - (
\nabla\Delta^S\psi,\nabla\psi) +|\nabla^2\psi|^2
\\&&+\sum_{\alpha,\beta=1}^n\Big(2(R^S(s_\alpha,s_\beta) \nabla_{s_\alpha}\psi,
\nabla_{s_\beta}\psi) +( ( \nabla_{s_\alpha}
R^S)(s_\alpha,s_\beta)\psi,\nabla_{s_\beta}\psi) \\&&\qquad\qquad+
\Ric(s_\alpha,s_\beta)(\nabla_{s_\alpha}\psi,\nabla_{s_\beta}\psi)\Big)\;.
\end{eqnarray*}
\end{Lemma}
\Proof With respect to $s_1,\dots,s_n$, the vector field $Y$ is
\[Y=\frac 12\sum_{\alpha, \beta=1}^n \nabla_{s_\alpha}\langle
\nabla_{s_\beta}\psi,\nabla_{s_\beta}\psi\rangle s_\alpha=\sum_{\alpha,
\beta=1}^n (\nabla^2_{s_\alpha,s_\beta}\psi,\nabla_{s_\beta}\psi) s_\alpha\,.\]
Let $x$ be a point in $M$.  For simplicity we choose $s_1,\dots,s_n$ such that
$\nabla s_\alpha=0$ in $x$ for all $\alpha =1,\dots,n$.  Then at $x$,
\begin{eqnarray*} \Div Y &=& \sum _{\alpha,\beta}^n
s_\alpha(\nabla^2_{s_\alpha,s_\beta}\psi,\nabla_{s_\beta}\psi) \\ &=&
\sum_{\alpha,\beta}^n \Big
((\nabla_{s_\alpha}\nabla_{s_\alpha}\nabla_{s_\beta}\psi, \nabla_{s_\beta}\psi)
- (\nabla_{s_\alpha}\nabla_{\nabla_{s_\alpha}s_\beta}\psi, \nabla_{s_\beta}
\psi) + (\nabla_{s_\alpha}\nabla_{s_\beta}\psi,
\nabla_{s_\alpha}\nabla_{s_\beta}\psi)\Big)\\
&=&\sum_{\alpha,\beta}^n\Big((\nabla_{s_\alpha}(\nabla_{s_\beta}
\nabla_{s_\alpha}\psi +\nabla_{[s_\alpha,s_\beta]}\psi +
R^S(s_\alpha,s_\beta)\psi) ,\nabla_{s_\beta}\psi)\\ &&\qquad
-(\nabla_{s_\alpha}\nabla_{\nabla_{s_\alpha}s_\beta}\psi, \nabla_{s_\beta}\psi) +
(\nabla_{s_\alpha}\nabla_{s_\beta}\psi, \nabla_{s_\alpha}\nabla_{s_\beta}\psi)
\Big)\\ &=& \sum_{\alpha,\beta}^n\Big(( \nabla_{s_\beta}\nabla_{s_\alpha}
\nabla_{s_\alpha}\psi + R^S(s_\alpha,s_\beta)\nabla_{s_\alpha}\psi +
\nabla_{s_\alpha}\nabla_{[s_\alpha,s_\beta]}\psi +
(\nabla_{s_\alpha}R^S)(s_\alpha,s_\beta)\psi\\ &&\qquad
+R^S(s_\alpha,s_\beta)\nabla_{s_\alpha}\psi-\nabla_{s_\alpha}
\nabla_{\nabla_{s_\alpha}s_\beta}\psi,\nabla_{s_\beta}\psi) + (
\nabla^2_{s_\alpha,s_\beta}\psi, \nabla^2_{s_\alpha,s_\beta}\psi)\Big)\\
&=&\sum_{\beta=1}^n(-\nabla_{s_\beta}\Delta^S\psi, \nabla_{s_\beta}\psi)
+\sum_{\alpha,\beta}^n\Big((-\nabla_{s_\beta}(\Div s_\alpha \cdot
\nabla_{s_\alpha}\psi)+ 2R^S(s_\alpha,s_\beta)\nabla_{s_\alpha}\psi \\
&&\qquad-
\nabla_{s_\alpha}\nabla_{\nabla_{s_\beta}s_\alpha}\psi +
(\nabla_{s_\alpha}R^S)(s_\alpha,s_\beta)\psi,\nabla_{s_\beta}\psi)+(
\nabla^2_{s_\alpha,s_\beta}\psi, \nabla^2_{s_\alpha,s_\beta}\psi)\Big)\,.
\end{eqnarray*}
Therefore it remains to show that at $x$,
\begin{equation}\label{ric}
\sum_{\alpha=1}^n\Big( -\nabla_{s_\beta}(\Div s_\alpha\cdot
\nabla_{s_\alpha}\psi)-\nabla_{s_\alpha}\nabla_{\nabla_{s_\beta}
s_\alpha}\psi\Big)=\sum_{\alpha=1}^n\Ric(s_\alpha,s_\beta)\nabla_{s_\alpha}
\psi\;.
\end{equation}
Since the $s_\alpha$ are orthonormal,
\[-s_\beta(\Div s_\alpha)=-s_\beta(\sum_{\gamma=1}^n
g(\nabla_{s_\gamma}s_\alpha,s_\gamma))=s_\beta(\sum_{\gamma=1}^n
g(s_\alpha,\nabla_{s_\gamma}s_\gamma))=\sum_{\gamma=1}^n
g(s_\alpha,\nabla_{s_\beta}\nabla_{s_\gamma}s_\gamma) \;, \]
and thus \[-\sum_{\alpha=1}^n \nabla_{s_\beta}(\Div s_\alpha\cdot
\nabla_{s_\alpha}\psi)=\sum_{\alpha,\gamma=1}^n
g(s_\alpha,\nabla_{s_\beta}\nabla_{s_\gamma}s_\gamma)
\nabla_{s_\alpha}\psi=\sum_{\beta,\gamma=1}^n \nabla_{\nabla_{s_\beta}
\nabla_{s_\gamma}s_\gamma}\psi\;. \]
On the other hand, using that $\nabla_{s_\beta} s_\alpha(x)=0$,
\begin{eqnarray*}
\nabla_{s_\alpha}\nabla_{\nabla_{s_\beta}s_\alpha}\psi
&=& \nabla_{\nabla_{s_\beta}s_\alpha} \nabla_{s_\alpha} \psi
\:+\: R^S(s_\alpha, \nabla_{s_\beta} s_\alpha) \:\psi \:+\:
\nabla_{[s_\alpha, \nabla_{s_\beta}s_\alpha]}\psi \\
&=& \nabla_{[s_\alpha, \nabla_{s_\beta}s_\alpha]}\psi
\;=\; \nabla_{\nabla_{s_\alpha}\nabla_{s_\beta} s_\alpha}\psi \;.
\end{eqnarray*}
Hence at $x$,
\begin{eqnarray*}
\lefteqn{ \sum_{\alpha=1}^n\Big(
-\nabla_{s_\beta}(\Div s_\alpha\cdot
\nabla_{s_\alpha}\psi)-\nabla_{s_\alpha}\nabla_{\nabla_{s_\beta}
s_\alpha}\psi\Big) } \\
&=& \sum_{\alpha=1}^n(\nabla_{\nabla_{s_\beta}\nabla_{s_\alpha}
s_\alpha}\psi - \nabla_{\nabla_{s_\alpha}\nabla_{s_\beta}
s_\alpha}\psi) \;=\; \sum_{\alpha=1}^n\nabla_{R(s_\beta,s_\alpha)s_\alpha}\psi
=\sum_{\alpha=1}^n\Ric(s_\alpha,s_\beta)\nabla_{s_\alpha}\psi \;.\;\;\;
\hspace*{.01cm} {\mbox{ \hfill}} \FBox
\end{eqnarray*}

\begin{Corollary} \label{coroll42}
For $\psi$ a solution of (\ref{eq:1C}) with $|\psi_0|=1$ and $\eta$ a positive
smooth function with $\sup_M (|\eta| + |\nabla \eta|)<\infty$,
\[ \int_M \eta \:|\nabla^2\psi|^2\,d\mu\;\leq\;
m \:C_1(n)\: \sup_M(|\eta R| + |\Delta \eta|) \:+\:
\sqrt{m}\:C_2(n)\: \|\eta \:\nabla R\|_2 \;. \]
\end{Corollary}
\Proof
We multiply the result of Lemma~\ref{L41} by $\eta$ and integrate over the
ball $D_\rho$, (\ref{eq:Dr}). Using Gauss' theorem and the definition of $Y$,
\begin{eqnarray}\label{int}
\lefteqn{\int_{D_\rho}\eta|\nabla^2\psi|^2\,d\mu \;=\; \int_{S_\rho} \left( \eta \: g(Y,\nu) -\frac 12 |\nabla \psi|^2 \: g(\nabla \eta,\nu)
\right) dS_\rho} \nonumber \\
&&-\frac 12 \int_{D_\rho} |\nabla \psi|^2\cdot\Delta \eta \,d\mu \:+\:
\int_{D_\rho} \eta(\nabla\Delta^S\psi,\nabla\psi) d\mu \nonumber\\
&& -\int_{D_\rho} \eta\sum_{\alpha,\beta=1}^n\Big(2(R^S(s_\alpha,s_\beta) \nabla_{s_\alpha}\psi,
\nabla_{s_\beta}\psi) +( ( \nabla_{s_\alpha}
R^S)(s_\alpha,s_\beta)\psi,\nabla_{s_\beta}\psi) \nonumber\\&&
\qquad\qquad\spc+
\Ric(s_\alpha,s_\beta)(\nabla_{s_\alpha}\psi,\nabla_{s_\beta}\psi)\Big)
\:d\mu \;.
\end{eqnarray}
Taking the limit $\rho \rightarrow\infty$, the integral over $S_\rho$ tends
to zero because of the asymptotic behavior of $\psi$ (\ref{psi}).
Furthermore, we can estimate the remaining terms in (\ref{int}) according to
\begin{eqnarray*}
\sum_{\alpha,\beta=1}^n|2(R^S(s_\alpha,s_\beta)\nabla_{s_\alpha}\psi,
\nabla_{s_\beta}\psi)|&=&\frac 12\sum_{\alpha,\beta,\gamma,\delta=1}^n|
R(s_\alpha,s_\beta,s_\gamma,s_\delta)s_\gamma\cdot s_\delta\cdot
\nabla_{s_\alpha}\psi,\nabla_{s_\beta}\psi)|\\
&\le&\tilde c_1(n)|R|\,|\nabla\psi|^2\\[2ex]
\sum_{\alpha,\beta=1}^n|((\nabla_{s_\alpha}R^S)(s_\alpha,s_\beta)\psi,
\nabla_{s_\beta}\psi)|&=&\frac14\sum_{\alpha,\beta,\gamma,\delta=1}^n|
((\nabla_{s_\alpha}R)(s_\alpha,s_\beta,s_\gamma,s_\delta)s_\gamma\cdot s_\delta
\cdot \psi,\nabla_{s_\beta}\psi)|\\
&\le&\tilde c_2(n) |\nabla R|\,|\psi|\,|\nabla \psi|\\[2ex]
\sum_{\alpha,\beta=1}^n \Ric(s_\alpha,s_\beta)(\nabla_{s_\alpha}\psi,
\nabla_{s_\beta}\psi)&\le&\tilde c_3(n)|R|\,|\nabla \psi|^2\\[2ex]
|(\nabla\Delta^S\psi,\nabla\psi)| &=&\sum_{\alpha=1}^n (\nabla_{s_\alpha}
(\frac \tau 4\,\psi),\nabla_{s_\alpha}\psi)|\\
&\le&\tilde c_4(n)|\nabla R|\, |\psi|\,|\nabla\psi|+\tilde c_5(n)|R|\,|\nabla\psi|^2.
\end{eqnarray*}
with suitable constants $\tilde c_1,\ldots,\tilde c_4$. We conclude that
\begin{eqnarray*}
\int_{D_r}\eta|\nabla^2\psi|^2\,d{D_r}&\le&\frac 12\int_{D_r}
|\nabla^2\psi|^2|\Delta\eta|\,d{D_r}\\&& + \tilde C_1(n) \int_{D_r} |\eta R|
\,|\nabla\psi|^2\,d{D_r}
 + \tilde C_2(n) \int_{D_r}\eta|\nabla R|\,|\psi|\,|\nabla \psi| \,d{D_r}\\
 &\le& \frac 12 \|\nabla\psi\|_2^2\, \sup_M|\Delta\eta|\\&&+\tilde C_1(n)
 \|\nabla\psi\|_2^2\, \sup_M |\eta R| + \tilde C_2(n) \|\eta\nabla R\|_2\,
 \|\nabla\psi\|_2\,\sup_M|\psi| \;,
\end{eqnarray*}
and the assertion follows from (\ref{eq:1C2}) and (\ref{eq:1E}).
\QED

\section{A-Priori Estimates for the Spinor Operator}
We choose an orthonormal basis of constant spinors $(\psi^i_0)_{i=1,\ldots,
N}$, $N=2^{[n/2]}$, $\langle \psi_0^i, \psi_0^j \rangle = \delta^{ij}$,
and denote the corresponding solutions of the boundary problem~(\ref{eq:1C})
by $(\psi^i)_{i=1,\ldots,N}$. We define the {\em{spinor operator}} $P_x$ by
\begin{equation}
    P_x \;:\; S_xM \longrightarrow S_xM \;:\; \psi \longmapsto
    \sum_{i=1}^N \langle \psi^i_x,\: \psi \rangle\: \psi^i_x \;.
    \label{eq:sop}
\end{equation}
Since at infinity the $\psi^i$ go over to an orthonormal basis,
\[ \lim_{|x| \rightarrow \infty} P_x \;=\; \1\;. \]
This section is concerned with a-priori estimates for the operator $P_x$.
In the following lemma, we use the maximum principle to derive upper bounds
for $|P_x|$, where $|\:.\:|$ denotes the sup-norm.
\begin{Lemma} \label{lemma31} $|P_x| \le 1$ for all $x\in M$.
\end{Lemma}
\Proof
We define the matrix $H$ by
\[ H \;=\; (h_{ij})_{i=1,\dots, N\,j=1,\dots,N} \spc{\mbox{with}}\spc
h_{ij}=\langle \psi^i_x,\psi^j_x\rangle\;. \]
By definition, $H$ is Hermitian and all eigenvalues of $H$
are real and non-negative. Consider the spinor $\psi=\sum_{i=1}^Na_i\psi^i$
with $a=(a_1,\dots,a_N)^{\top}\in \C^N$. Then $\psi$ is a solution of the
boundary problem~(\ref{eq:1C}) with $\psi_0=\sum_{i=1}^Na_i\psi^i_0$.
By (\ref{eq:1E}),
\[ \sum_{i,j=1}^N\overline{a_i}\:a_j\:h_{ij} \;=\;
|\psi_x|^2 \;\leq\; |\psi_0|^2 \;=\; \sum_{i=1}^N |a_i|^2 \]
for all $x \in M$.
Since $a$ is arbitrary, all eigenvalues of $H$ must be $\le 1$.
Let $\phi$ be an arbitrary spinor at $x$ and $\psi$ its orthogonal
projection to $\Span\{\psi^i\mid i=1,\dots,N\}$.
Then $\psi=\sum_{i=1}^N a_i \psi^i$ and
\begin{eqnarray*}
|P_x \phi|^2 &=& \sum_{i,j=1}^N \langle \psi,\psi^i_x\rangle\:
\langle\psi^i_x,\psi^j_x\rangle \: \langle\psi^j_x, \psi \rangle
\;=\;\overline{(H a)}^{\top}H (H a) \\
&\le& \overline{a}^{\top}Ha \;=\; |\psi_x|^2 \;=\; |\phi_x|^2 \;,
\end{eqnarray*}
where the first inequality uses that there is an orthonormal basis of $\C^N$ of
eigenvectors of $H$ and that all eigenvalues of $H$ are non-negative and
$\le 1$.
\QED
We next derive Sobolev estimates for the Hilbert--Schmidt norm $\| \cdot \|$ of the
operator $\1-P_x$.
\begin{Lemma} \label{lemma32}
There is a constant $c$ depending only on the dimension $n$ such that
for any $\varepsilon>0$,
\begin{equation} \label{eq:3A}
\|\1 - P_x \|^2 \;<\; \varepsilon
\end{equation}
except on a set $D(\varepsilon)$ with
\begin{equation} \label{eq:3B}
\mu(D) \;\leq\; \left( \frac{c\:m}{\varepsilon^2\:k^2}
\right)^{\frac{n}{n-2}} \;\;\;.
\end{equation}
\end{Lemma}
{\Proof}
We set $h(x) = \|\1-P_x\|^2$. Choosing an orthonormal basis
$(\phi^j)_{j=1,\ldots,N}$
of $S_xM$, the trace of $S_x$ is computed as follows,
\[ \Tr (P_x) \;=\; \sum_{j=1}^N \langle \phi^j, P_x \:\phi^j \rangle
\;\stackrel{(\ref{eq:sop})}{=}\; \sum_{i,j=1}^N \langle \psi^i, \phi^j \rangle
\langle \phi^j, \psi^i \rangle \;=\; \sum_{i=1}^N \langle \psi^i,\psi^i \rangle
\;. \]
Similarly we obtain for $h$,
\[ h(x) \;=\; \Tr (\1 - 2 P_x + P_x^2) \;=\;
N \:-\: 2 \sum_{i=1}^N \langle \psi^i_x,\:\psi^i_x \rangle \:+\:
\sum_{i,j=1}^N |\langle \psi^i_x,\:\psi^j_x \rangle|^2 \;, \]
and the gradient of $h$ is computed to be
\[ \nabla h \;=\; -4\: \sum_{i=1}^N  \sum_{\alpha=1}^n \langle
\nabla_{s_\alpha} \psi^i,\:
\psi^i \rangle \:s_\alpha \:+\: 4\: \sum_{i=1}^N\sum_{\alpha=1}^n
\langle \nabla_{s_\alpha} \psi^i,\: P_x\: \psi^i \rangle \: s_\alpha \;. \]
The Schwarz inequality combined with Lemma~\ref{lemma31} yield that
\[ |\nabla h| \;\leq\; 8 n \sum_{i=1}^N |\nabla \psi^i| \;, \]
and applying H\"older's inequality,
\[ |\nabla h|^2 \;\leq\; 64 n^2\:N \:\sum_{i=1}^N |\nabla \psi^i|^2 \;. \]
We now integrate over $M$ and substitute in~(\ref{eq:1C2}) to obtain
\[ \|\nabla h\|_2^2 \;\leq\; 16 n^2\:N^2\: c(n)\; m \;. \]
The Sobolev inequality yields for some constant $C(n)$,
\begin{equation}
    k^2\:\| h \|_q^2 \;\leq\; C\:\| \nabla h \|_2^2 \spc {\mbox{with}}\spc
    q \;=\; \frac{2n}{n-2}\;, \label{eq:4sobolev}
\end{equation}
where $k$ is the isoperimetric constant of $M$ (see Lemma~\ref{lemma_sob} for
the derivation of this inequality).
Hence $h<\varepsilon$
except on a set $D$ of small measure~(\ref{eq:3B}).
\QED

For the reader not familiar with Sobolev inequalitities in non-compact
Riemannian manifolds we now give the proof of
inequality~(\ref{eq:4sobolev}).
\begin{Lemma} \label{lemma_sob}
Let $0 \leq h \in C^\infty(M)$ with $\lim_{|x| \rightarrow \infty} h(x)=0$
and $q=2n/(n-2)$. Then
\[ \| h \|_q \;\leq\; \frac{q}{k}\: \|\nabla h \|_2\;. \]
\end{Lemma}
{\Proof}
We define for $u>0$ the sets
\[ N_u \;=\; \{ x \in M \:|\: h(x) \geq u \} \;\;\;,\spc S_u \;=\; \partial
N_u \;=\; h^{-1}(u)\;. \]
Since $\lim_{|x| \rightarrow \infty} h(x)=0$, these sets are compact.
Sard's lemma yields that, with the exception of $u$ in a set of measure
zero, $\nabla h$ does not vanish on $S_u$, and so $S_u$ is a compact
submanifold of $M$ of codimension one. We denote the volume of $N_u$ by
$V_u$ and the area of $S_u$ by $A_u$. Also, $dS_u$ denotes the measure
on $S_u$ corresponding to the induced Riemannian metric.
The co-area formula yields that for any $p>0$,
\begin{equation}
\int_M h^p \:|\nabla h| \:d\mu \;=\;
\int_0^\infty du \:\int_{S_u} h^p \:dS_u \;=\;
\int_0^\infty u^p\:A_u\:du\;. \label{eq:4trans}
\end{equation}
Furthermore, we have the following estimates,
\begin{eqnarray}
V_u &=& \int_{N_u} d\mu \;\leq\; \frac{1}{u^q} \:\int_{N_u} h^q\: d\mu
\;\leq\; u^{-q} \:\|h\|_q^q \label{eq:432} \\
\|h\|_q^q &=& \int_M h^q \:d\mu \;=\; q \int_M \left( \int_0^{h(x)} u^{q-1}
\:du \right) d\mu_x \nonumber \\
&=& q \int_M d\mu_x \int_0^\infty du\: u^{q-1}\: \Theta(h(x)-u) \;,
\label{eq:439}
\end{eqnarray}
where $\Theta$ is the Heaviside function $\Theta(x)=1$ for $x \geq 0$ and
$\Theta(x)=0$ otherwise. The integrand in~(\ref{eq:439}) is positive,
and thus we may commute the integrals according to Fubini's theorem,
\begin{eqnarray*}
\|h\|_q^q &=& q \int_0^\infty du\: u^{q-1} \int_M \Theta(h(x)-u)\: d\mu_x \\
&=& q \int_0^\infty u^{q-1} \:V_u \:du
\;=\; q \int_0^\infty u^{q-1} \:V_u^{\frac{1}{n}}\: V_u^{\frac{n-1}{n}} \:du
\;.
\end{eqnarray*}
The isoperimetric inequality bounds the factor $V_u^{\frac{n-1}{n}}$ from
above by $A_u/k$. The factor $V^{\frac{1}{n}}$, on the other hand, can be
estimated with~(\ref{eq:432}). We thus obtain the inequality
\[ \|h\|_q^q \;\leq\; \frac{q}{k}\: \|h\|_q^{\frac{q}{n}} \:\int_0^\infty
u^{\frac{q}{2}}\:A_u\: du \;\;\;. \]
We finally substitute in~(\ref{eq:4trans})
and apply the Schwarz inequality,
\[ \;\;\;\;\;\;\;\;
\spc\spc \|h\|_q^q \;\leq\; \frac{q}{k}\: \|h\|_q^{\frac{q}{n}} \:\int_M h^{\frac{q}{2}}
\:|\nabla h| \:d\mu \;\leq\; \frac{q}{k}\: \|h\|_q^{q-1}\; \|\nabla h\|_2
\;\;\;. \spc\spc {\mbox{ \hfill}} \FBox \]

\section{Proof of the Curvature Estimates}
We derive a pointwise estimate of the curvature tensor in terms of
the spinors $\psi^i$ and their second derivatives.
\begin{Lemma} \label{lemma51}
\[ \left( N \:-\: \sqrt{8N} \:\| \1-P_x \| \right)
|R|^2 \;\leq\; 32 \sum_{i=1}^N |\nabla^2 \psi^i|^2 \;. \]
\end{Lemma}
{\Proof}
For convenience we again choose an orthonormal frame $s_1,\ldots,s_n$ with
$\nabla s_\alpha(x)=0$. The definition of curvature~(\ref{eq:2Z}) and the
Schwarz inequality give the following relation between curvature and the
second derivatives of $\psi$,
\begin{eqnarray}
\lefteqn{ -\sum_{\alpha, \beta=1}^n \langle \psi,\:
R^S(s_\alpha, s_\beta)^2\:\psi \rangle \;=\; \sum_{\alpha,\beta=1}^n
\langle R^S(s_\alpha,s_\beta)\psi,R^S(s_\alpha,s_\beta)\psi\rangle }
\nonumber \\
&=&\sum_{\alpha,\beta=1}^n \langle
\nabla^2_{s_\alpha,s_\beta}\psi-\nabla^2_{s_\beta,s_\alpha}\psi,\nabla^2_{s_\alpha,
s_\beta}\psi-\nabla^2_{s_\beta,s_\alpha}\psi\rangle \nonumber \\
&=& \sum_{\alpha,\beta=1}^n (|\nabla^2_{s_\alpha,s_\beta}\psi|^2 +
|\nabla^2_{s_\beta,s_\alpha}\psi|^2-2\langle
\nabla^2_{s_\alpha,s_\beta}\psi,\nabla^2_{s_\beta,s_\alpha}\psi\rangle
\nonumber \\
&\le& 2\sum_{\alpha,\beta=1}^n (|\nabla^2_{s_\alpha,s_\beta}\psi|^2 +
|\nabla^2_{s_\beta,s_\alpha}\psi|^2)\ =\ 4|\nabla^2\psi| \;.
\label{eq:5b}
\end{eqnarray}
Using~(\ref{eq:2Y}) and~(\ref{eq:acomm}), the square of spinor curvature
is computed to be
\begin{eqnarray}
\lefteqn{ \sum_{\alpha, \beta=1}^n R^S(s_\alpha, s_\beta)^2
\;=\; \frac{1}{16} \sum_{\alpha, \beta=1}^n \;\sum_{\gamma,\delta,\epsilon,
\rho=1}^n R(s_\gamma, s_\delta, s_\alpha, s_\beta)\:
R(s_\epsilon, s_\rho, s_\alpha, s_\beta)\:
s_\gamma s_\delta s_\epsilon s_\rho } \\
&=& -\frac{1}{8} \:|R|^2 \:+\:
\frac{1}{16} \sum_{\alpha, \beta=1}^n \;\sum_{\gamma,\delta,\epsilon,
\rho {\mbox{\scriptsize{ all different}}}}
 R(s_\gamma, s_\delta, s_\alpha, s_\beta)\:
R(s_\epsilon, s_\rho, s_\alpha, s_\beta)\:
s_\gamma s_\delta s_\epsilon s_\rho \;. \spc \label{eq:5a}
\end{eqnarray}
In dimension $n \leq 3$, the second term in~(\ref{eq:5a}) clearly vanishes,
and so~(\ref{eq:5a}) is a multiple of the identity matrix, making it possible
to proceed as in~\cite{BF}. In order to control the second term
in~(\ref{eq:5a}), we consider the expectation value with respect to all
$\psi^i$s and take their sum,
\begin{eqnarray}
\lefteqn{ -\sum_{i=1}^N \sum_{\alpha, \beta=1}^n \langle \psi^i,\:
R^S(s_\alpha, s_\beta)^2 \:\psi^i \rangle \;=\;
-\sum_{\alpha, \beta=1}^n \Tr \left(R^S(s_\alpha, s_\beta)^2\:P_x \right) }
\nonumber \\
&=&-\sum_{\alpha, \beta=1}^n \Tr \left(R^S(s_\alpha, s_\beta)^2 \right)
\:+\: \sum_{\alpha, \beta=1}^n \Tr
\left(R^S(s_\alpha, s_\beta)^2\:(\1 - P_x) \right) \nonumber \\
&\geq& -\sum_{\alpha, \beta=1}^n \Tr \left(R^S(s_\alpha, s_\beta)^2 \right)
\:-\: \|\sum_{\alpha, \beta=1}^n R^S(s_\alpha, s_\beta)^2 \|\;
\|\1-P_x\| \;. \label{eq:5c}
\end{eqnarray}
A straightforward calculation using~(\ref{eq:5a}) shows that
\[ -\sum_{\alpha, \beta=1}^n \Tr \left(R^S(s_\alpha, s_\beta)^2 \right)
\;=\; \frac{N}{8} \:|R|^2 \;\;\;,\spc
\|\sum_{\alpha, \beta=1}^n R^S(s_\alpha, s_\beta)^2 \|
\;\leq\; \sqrt{\frac{N}{8}}\:|R|^2 \;. \]
Substituting these formulas into~(\ref{eq:5c}) and using~(\ref{eq:5b}) gives
the result.
\QED
Our main result follows by combining Lemma~\ref{lemma32},
Corollary~\ref{coroll42}, and Lemma~\ref{lemma51}.\\[.5em]
{\em{Proof of Theorem~\ref{thm1}}:}
We choose the set $D$ as in Lemma~\ref{lemma32} corresponding to
$\varepsilon=N/32$. Then according to~(\ref{eq:3A}) and Lemma~\ref{lemma51},
\[ \int_{M \setminus D} \eta\:|R|^2\:d\mu \;\leq\;
\int_{M \setminus D} \eta \;\frac{64}{N}\:\sum_{i=1}^N
|\nabla^2 \psi^i|^2\:d\mu \;\leq\;
\frac{64}{N}\:\sum_{i=1}^N \int_M \eta \:|\nabla^2 \psi^i|^2\:d\mu \;, \]
where in the last step we used the positivity of the integrand. Now
apply Corollary~\ref{coroll42}.
\hspace*{1.5cm} \hfill \QED

\noindent{\em{Acknowledgments:}} We would like to thank Robert Bartnik and
Helga Baum for helpful suggestions and comments. We are grateful to Hubert
Bray for many inspiring discussions and valuable ideas.

Max Planck Institute for Mathematics in the Sciences,
Inselstr.\ 22-26, 04103 Leipzig, Germany, {\tt{Felix.Finster@mis.mpg.de}},
{\tt{ikath@mis.mpg.de}}


\addcontentsline{toc}{section}{References}
\begin{thebibliography}{MM}
\bibitem{SY} Schoen, R., Yau,  S.-T., ``On the proof of the positive
mass conjecture in General Relativity,'' {\em{Commun.\ Math.\ Phys.}}\ 65
(1976) 45-76
\bibitem{W} Witten, E., ``A new proof of the positive energy
theorem,'' {\em{Commun.\ Math.\ Phys.}}\ 80 (1981) 381-402
\bibitem{Br} Bartnik, R., ``The mass of an asymptotically flat manifold,''
{\em{Comm.\ Pure Appl.\ Math.}} 39 (1986) 661-693
\bibitem{BF} Bray, H., Finster, F., ``Curvature estimates and the
positive mass theorem,'' math.DG/9906047, to appear in
{\em{Comm.\ in Analysis and Geometry}} (2001)
\bibitem{PT} Parker, T., Taubes, C.\ H., ``On Witten's proof of the
positive energy theorem,'' {\em{Commun.\ Math.\ Phys.}}\ 84 (1982)
223-238
\bibitem{B} Baum, H., ``Spin-Strukturen und Dirac-Operatoren \"uber
pseudoriemannschen Mannigfaltigkeiten,'' {\em{Teubner Verlag Leipzig}} (1981)
\bibitem{LM} Lawson, H.B., Michelson M.-L., ``Spin Geometry,'' {\em{Princeton
 Univ.\ Press}} (1989)
\end{thebibliography}
\end{document}